# CORRECTION NOTE
## TYPICAL CONFIGURATION FOR ONE-DIMENSIONAL RANDOM FIELD KAC MODEL[1]

By Marzio Cassandro, Enza Orlandi and Pierre Picco

*Università di Roma, Università di Roma Tre and CPT-CNRS*

Estimate (3.39) which appears in the proof of Proposition 3.4 in [*Ann. Probab.* **27** (1999) 1414–1467] is wrong. We present below a corrected proof which introduces an extra factor 2 in equations (3.34) and (3.35). This has no consequence in the rest of the paper since Proposition 3.4 is used to estimate only ratios; see (3.23) and (3.25).

In Proposition 3.4 in [1], the condition $m \in \{-1, -1 + 2/|B|, -1 + 4/|B|, \ldots, 1 - 2/|B|, 1\}$ has to be added. This is harmless since Proposition 3.4 is used for proving Proposition 3.1, where this assumption is done. Moreover, (3.34) and (3.35) must be replaced respectively by

$$\Psi_{z,\alpha,m} = \frac{2}{\sqrt{2\pi|B|}\sigma_z}\left(1 \pm \frac{66}{|B|\sigma_z^2}\right) \tag{1.1}$$

and

$$\Psi_{z,\alpha,m} = \frac{2}{\sqrt{2\pi|B|}\sigma_z}\left(1 \pm \frac{66}{g(|B|)}\right). \tag{1.2}$$

Below we outline the arguments to get (1.2), the case of (1.1) is similar.

In the proof of Proposition 3.4, inequality (3.39) is clearly wrong for $k = \pm\pi$. Since, for $y \in [0,1]$, we have $|ye^{-2ik} + (1-y)|^2 = 1 - 2y(1-y)(1 - \cos(2k))$ and $1 - s \leq e^{-s}$ for all $s \in \mathbb{R}$, it is easy to see that

$$\left|\frac{\cosh(x \pm ik)}{\cosh(x)}\right| \leq \exp\left[-\frac{1 - \cos(2k)}{4\cosh^2 x}\right] \tag{1.3}$$

Received June 2005; revised September 2005.

[1]Supported by CNR-CNRS-Project 8.005, INFM-Roma; MURST/Cofin 01-02/03-04.

that replaces (3.39). Then, using $\cos(x) \leq 1 - \frac{x^2}{2} + \frac{x^4}{4!}$, it can be checked that, for $k \in [0, \pi]$,

$$(1.4) \qquad 1 - \cos(2k) \geq 2\left(1 - \frac{\pi^2}{12}\right)(k^2 \wedge (k - \pi)^2),$$

from which one gets, for $k \in [0, \pi]$,

$$(1.5) \qquad |\Phi(z, \alpha, k)| \leq \exp\left[-\frac{(1 - \pi^2/12)(k^2 \wedge (k - \pi)^2)}{2}|B|\sigma_z^2\right],$$

where $\Phi(z, \alpha, k)$ is defined in (3.38) and $\sigma_z$ is defined in (3.28) in [1]. Formula (1.5) replaces (3.40) in [1]. As a consequence, (3.41) has to be replaced by

$$(1.6) \qquad \widetilde{\mathcal{E}}_\rho = \frac{1}{2\pi}\int_{-\pi}^{\pi} \mathbb{1}_{\{\rho < |k| \leq \pi - \rho\}}\Phi(z, \alpha, -k)e^{ikm|B|}\, dk.$$

Then choosing as in [1], $\rho = (\sigma_z\sqrt{|B|})^{-1}f(|B|)$ with

$$(1.7) \qquad f(|B|) = \sqrt{\frac{2}{1 - \pi^2/12}\log g(|B|)},$$

where $g$ is as in Proposition 3.4 in [1], one gets

$$(1.8) \qquad |\widetilde{\mathcal{E}}_\rho| \leq \frac{1}{\sqrt{2\pi|B|}\sigma_z}\left(\frac{2}{\sqrt{\pi(1 - \pi^2/12)\log g(|B|)}}\right)\frac{1}{g(|B|)},$$

that replaces (3.48) in [1]. Calling as in [1] [see (3.45)],

$$(1.9) \qquad \Psi_{z,\alpha,m}(\rho) = \frac{1}{2\pi}\int_{-\rho}^{+\rho} e^{ik|B|m}\Phi(z, \alpha, k)\, dk,$$

introducing the two quantities

$$(1.10) \qquad \begin{aligned} I_2 &= \frac{1}{2\pi}\int_{-\pi}^{-\pi+\rho} e^{ik|B|m}\Phi(z, \alpha, k)\, dk, \\ I_3 &= \frac{1}{2\pi}\int_{\pi-\rho}^{\pi} e^{ik|B|m}\Phi(z, \alpha, k)\, dk. \end{aligned}$$

After simple algebra, using that $m = -1 + \frac{2l}{|B|}$ for some $l \in \mathbb{Z}$ and elementary change of variables, one gets the crucial relation

$$(1.11) \qquad I_2 + I_3 = \Psi_{z,\alpha,m}(\rho).$$

Now $\Psi_{z,\alpha,m}$ defined in (3.37) satisfies

$$(1.12) \qquad \Psi_{z,\alpha,m} = 2\Psi_{z,\alpha,m}(\rho) + \widetilde{\mathcal{E}}_\rho.$$

The extra factor 2 we mention in the abstract is the one in (1.12). Using the same computations done after (3.45) in [1], one gets (1.2).

M. Cassandro
Dipartimento di Fisica
Università di Roma "La Sapienza"
INFM-Sez. di Roma. P. le A. Moro
00185 Roma
Italy
E-mail: cassandro@roma1.infn.it

E. Orlandi
Dipartimento di Matematica
Università di Roma Tre
L.go S.Murialdo 1
00156 Roma
Italy
E-mail: orlandi@mat.uniroma3.it

P. Picco
CPT, UMR CNRS 6207
Université de Provence Aix–Marseille 1
Université de la Mediterranée Aix–Marseille 2
et Université de Toulon et du Var
Luminy, Case 907, 13288
Marseille Cedex 9
France
E-mail: picco@cpt.univ-mrs.fr